%% file: ttsurvey.tex
\documentclass[11pt,a4paper,oneside]{amsart}
\usepackage{amssymb,amsmath,amsfonts,latexsym,graphics,graphicx,color,vmargin}

\newtheorem{theorem}{Theorem}[section]

\newtheorem{conjecture}[theorem]{Conjecture}

\theoremstyle{definition}

\newcommand{\beql}[1]{\begin{equation}\label{#1}}
\newcommand{\eeq}{\end{equation}}
\newcommand{\comment}[1]{}

\newcommand{\Abs}[1]{{\left|{#1}\right|}}

\newcommand{\Set}[1]{{\left\{{#1}\right\}}}

\newcommand{\RR}{{\mathbb R}}
\newcommand{\CC}{{\mathbb C}}
\newcommand{\ZZ}{{\mathbb Z}}

\newcommand{\ft}[1]{\widehat{#1}}

\newcommand{\supp}{{\rm supp\,}}

\newcounter{rem}
\newcounter{step}
\begin{document}

\title{Tilings by translation}

\author{Mihail N. Kolountzakis \& M\'at\'e Matolcsi}

\date{September 2010}

\address{M.K.: Department of Mathematics, University of Crete, Knossos Ave.,
GR-714 09, Iraklio, Greece} \email{kolount@gmail.com}

\address{M.M.: Alfr\'ed R\'enyi Institute of Mathematics,
Hungarian Academy of Sciences POB 127 H-1364 Budapest, Hungary.
\newline
 (also at BME Department of Analysis, Budapest,
H-1111, Egry J. u. 1)}\email{matomate@renyi.hu}

\thanks{MK \& MM: Supported by research grant No 2937 from the Univ.\ of Crete
and by a grant of the Universita Autonoma de Madrid.}
\thanks{MM: Supported by the ERC-AdG 228005, and OTKA Grants No. K77748, K81658.}


\maketitle
\tableofcontents



\section{Introduction}\label{sec:intro}

To tile means to cover a given part of space, without overlaps, using a small
number, often just one, of different types of objects. When we are covering the
floor of a room using identical rectangular tiles (in the ordinary sense of the
word) we put down the copies of the tiles in a regular way next to each other,
leaving no gaps.

Your floor can be as boring as the rectangular floor on the left or as interesting as the Escher lizard tiling on the right:\\
\begin{center}
\begin{tabular}{cc}
\resizebox{!}{5cm}{\includegraphics{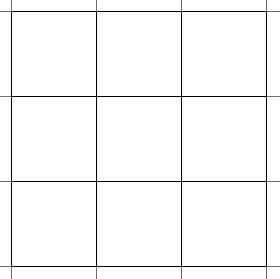}}
&
\resizebox{!}{5cm}{\includegraphics{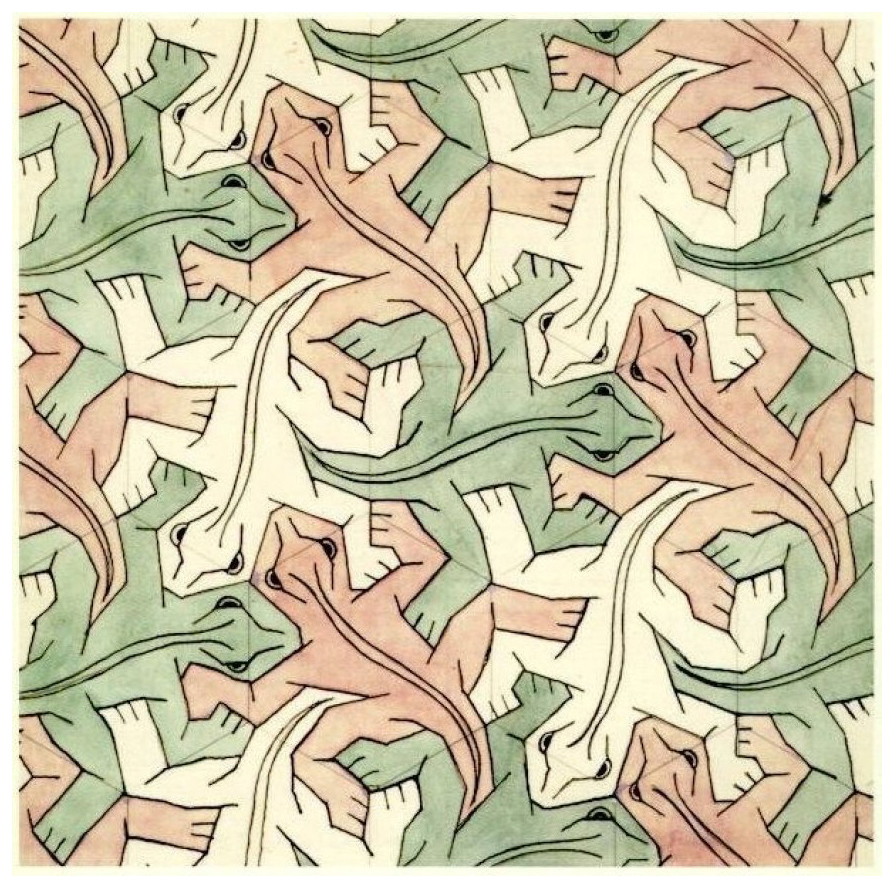}}
\end{tabular}
\end{center}
In both these floor tilings there is only one tile used, a square tile on the left and a lizard shape tile
\begin{center}
\begin{tabular}{cc}
\resizebox{!}{4cm}{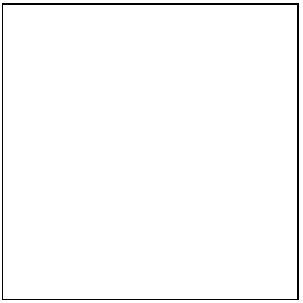}
&
\resizebox{!}{4cm}{\includegraphics{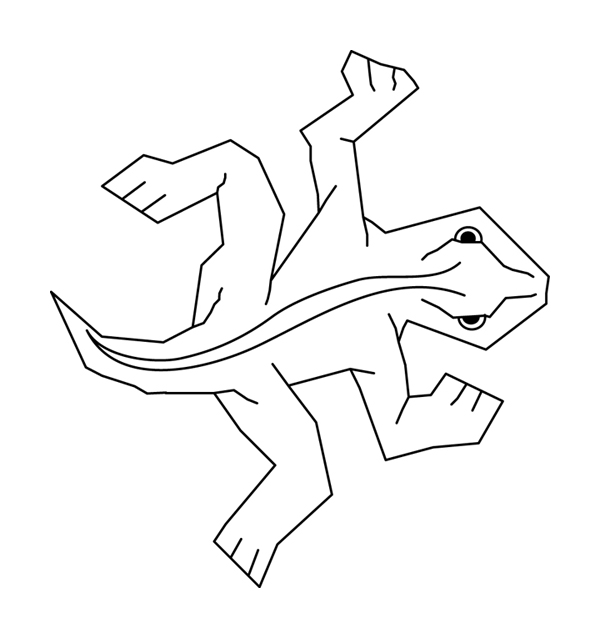}}
\end{tabular}
\end{center}
 on the right.

Beyond aesthetics and decoration of the tile there is another difference between these two examples
which is very significant from our point of view. To fill the floor using the square tile we only need to translate the tile around; we never have to turn it. This is not the case with the lizard tile. To fill the floor the various copies of the lizard tiles have to be turned appropriately so as to fit perfectly with each other.

It is with the seemingly boring case that we shall occupy ourselves in this survey: the {\em single} tile that we have at our disposal is {\em only allowed to translate} in space. We cannot rotate it or reflect it.
True, in the vast majority of examples in the literature (see e.g. \cite{grunbaum1986tilings}) that exhibit interesting behavior (such as undecidability or aperiodicity, on which we'll have more to say later) the allowed {\em group of motions} is usually the full group of rigid motions in space and therefore rotations are allowed (and some times even reflections as well). One would indeed be hard pressed to find eye-catching examples such as the Penrose ``kite and dart'' aperiodic set of tiles shown below.
\begin{center}
 \resizebox{!}{6cm}{\includegraphics{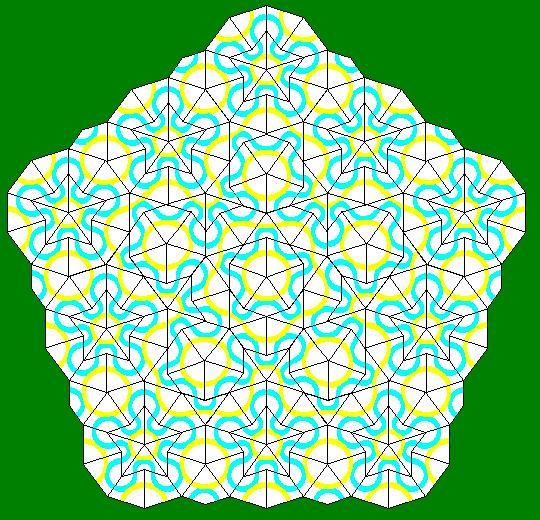}}
\end{center}
But we have good reasons for restricting ourselves to {\em tiling by
translation} which we hope to expose in this survey. To begin with let us state
that tiling by translation is still full of unresolved questions with
connections to number theory, Fourier analysis and the theory of computation,
and that happens even in dimension 1. If it is hard to imagine how one can
seriously study questions of tiling in dimension 1 let us point out that the
tiles need not have a ``nice'' shape. For instance, one can easily show by trial
and error that the set $E = \Set{0,2,3,5}$ (shown below in red)
\begin{center}
 \resizebox{4cm}{!}{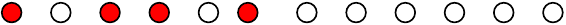}
\end{center}
can tile the integers $\ZZ$ by translation only and that such a tiling has
period 12. This trial and error process can become slower and slower if one
considers larger finite sets of integers and simply asks if they could be used
to tile $\ZZ$. In fact, it is not known how to do this in any
way that is essentially faster than trial and error.
This reflects how little we understand this tiling phenomenon.

Here is another example, of a more geometric flavor, to convince the reader that taking away the freedom to rotate and reflect the tile does not take away the fun. For this consider the so called {\em notched cube} (left)

\begin{center}
\begin{tabular}{cc}
 \resizebox{!}{4cm}{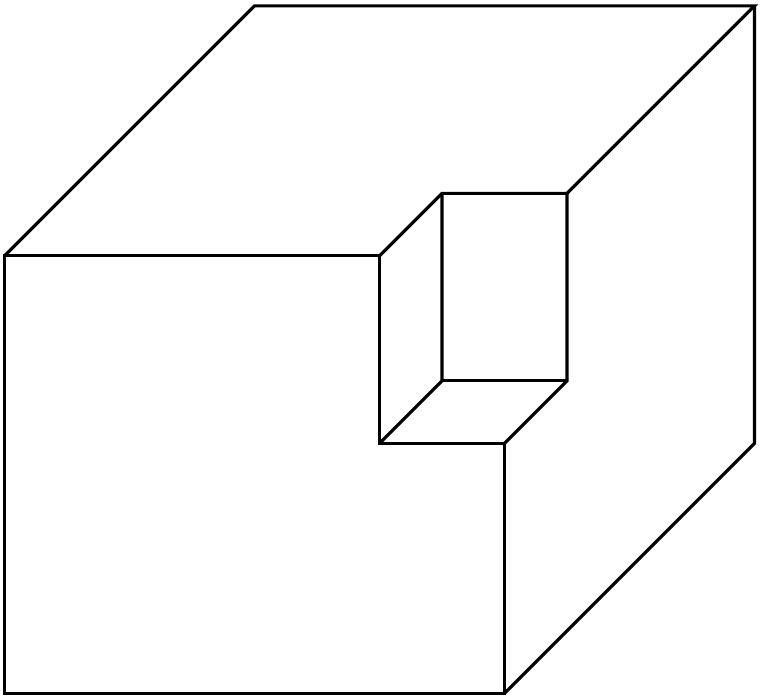}
&
 \resizebox{!}{4cm}{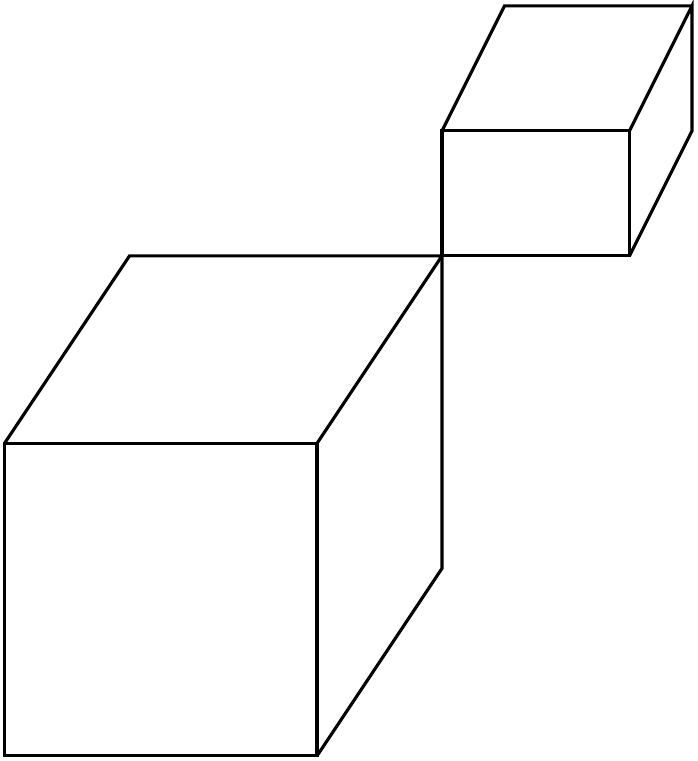}
\end{tabular}
\end{center}

\noindent
and try to visualize a way to fill 3 dimensional space with it (here it would not help to rotate or reflect it). We are making no assumptions here about the lengths of the rectangular cut (notch) at one corner of the
cube or rectangle. Whatever these lengths are this notched cube can indeed tile space by translation \cite{kolountzakis1998notched}. Yet the simplest way we know how to prove this is using Fourier analysis, making absolutely no use of geometric intuition.
Same holds for the {\em extended cube} (above right) and this seems to be even harder to visualize. At least the notched cube ``forces'' a corner of a copy of the tile to fill the notch in another copy and this can at least get you started but no such restriction is obvious for the extended cube.

By now we should fix the rules of the game. What does it mean for a subset
$\Omega$ of $\RR^d$ to tile $\RR^d$ when translated at the locations $\Lambda
\subseteq \RR^d$ (another {\em discrete} set)? Simply that the copies
$$
\Omega + \lambda,\ \ \lambda \in \Lambda,
$$
are mutually {\em non-overlapping}, and {\em cover the whole space}.
Here, as is customary, we denote by
$\Omega+\lambda$ the set
$$
\Set{\omega+\lambda:\ \omega \in \Omega},
$$
\begin{center}
 \def\pointslambda{$\Lambda$}
 \resizebox{6cm}{!}{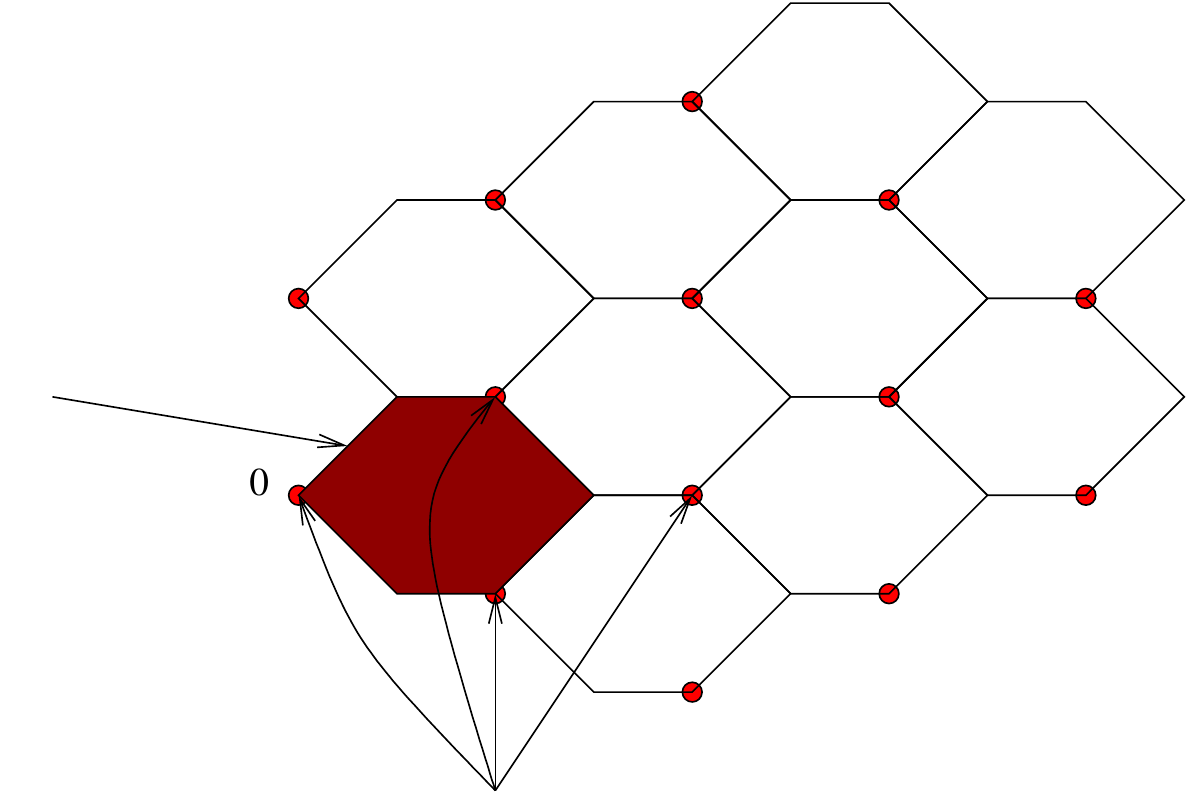}
\end{center}
or, in other words, the translate of $\Omega$ by the vector $\lambda$. And what
does it mean for two sets to be non-overlapping? One could, for instance, demand
that the interiors, but not the boundaries, of the two sets are disjoint. But
it turns out that the most convenient way is to demand that the intersection of
these two sets has zero volume (or area in dimension $d=2$, or length in
dimension $d=1$) or, more precisely, zero Lebesgue measure.

\section{Using the Fourier Transform}\label{sec:fourier}

\subsection{Tiling in the Fourier domain}\label{sec:tiling-in-fourier}

Having specified what we mean we turn now to exploiting the technology that we
know, and in this case we find it very fruitful to define tiling simply by
the equation
\beql{disjoint-translates}
\sum_{\lambda \in \Lambda} \chi_\Omega(x-\lambda) = 1,\ \
\mbox{(almost everywhere for $x\in \RR^d$)}.
\eeq
We now have an equation! as a physicist would exclaim. A mathematician should
never underestimate the power of formal manipulation and the phenomena it could
reveal. So, let us rewrite our equation (ignoring from now on the exception of $0$
measure):
\beql{convolution}
\chi_\Omega * \delta_\Lambda = 1.
\eeq
Here the $*$ operator denotes convolution and $\delta_\Lambda$ is the {\em
measure}

\begin{center}
\resizebox{6cm}{2cm}{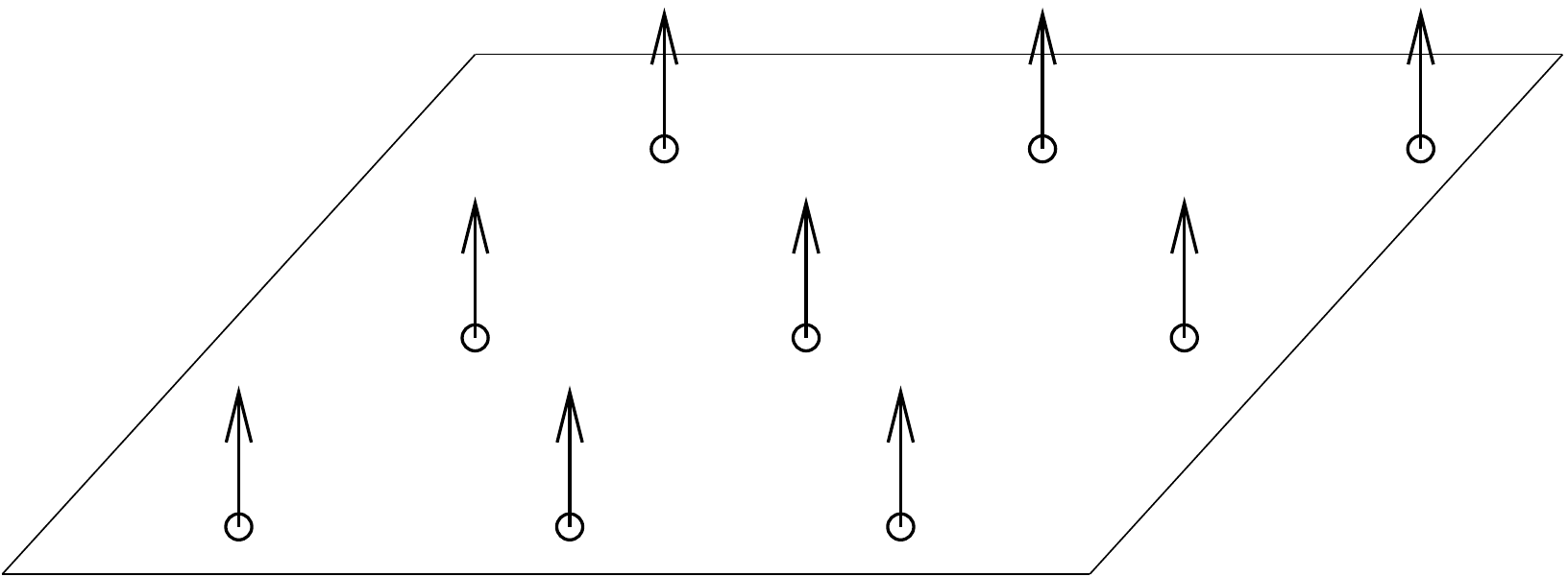}
\end{center}
$$
\delta_\Lambda = \sum_{\lambda \in \Lambda} \delta_\lambda,
$$
where $\delta_\lambda$ is a unit point mass sitting at point $\lambda \in
\RR^d$. (We recall that $f*\delta_\lambda (x) = f(x-\lambda)$, so that
convolving a function $f$ with $\delta_\lambda$ merely translates the function
by the vector $\lambda$.)
The object $\delta_\Lambda$ is a compact way to encode all information that's
contained in the set of translates $\Lambda$ in a way that we can operate
algebraically on it.

An equation such as \eqref{convolution} of course begs to be subjected to the
Fourier Transform\footnote{This of course depends on what training {\em you} have been
subjected to.} as the Fourier Transform (FT) behaves so nicely with convolution
(the FT of a convolution is the pointwise product of the FTs of the
convolution factors: $\ft{f*g} = \ft{f}\ft{g}$). The definition of the FT
$\ft{f}$ of a function $f$ that we use is
\beql{ft-def}
\ft{f}(\xi) = \int_{\RR^d} f(x) e^{-2\pi i \xi\cdot x}\,dx,\ \ (\xi\in\RR^d,
\int_{\RR^d}\Abs{f} < \infty).
\eeq
This definition is sufficient for the function $f=\chi_\Omega$ which is
integrable but one needs to define the FT of a measure by duality (see
for instance \cite{kolountzakis2004milano} for the details). Ignoring
such subtleties we take the FT of both sides of \eqref{convolution}
to get
\beql{conv-ft}
\ft{\chi_\Omega} \cdot \ft{\delta_\Lambda} = \delta_0.
\eeq
The measure $\delta_0$ on the right is just a point mass at point $0$. The
support (where it is ``nonzero'') of this is $\Set{0}$ and therefore that must
be the support of the left hand side of \eqref{conv-ft}. Since that is a product
it follows that, apart from point $0$, wherever $\ft{\chi_\Omega}$ is nonzero,
$\ft{\delta_\Lambda}$ must be zero so as to kill the product. Let us summarize
this in the following necessary (and often sufficient) condition for tiling:
\beql{ft-cond}
\supp \ft{\delta_\Lambda} \subseteq \Set{\ft{\chi_\Omega}=0} \cup \Set{0}.
\eeq
A formal proof of \eqref{ft-cond} can be looked up in
\cite{kolountzakis2004milano}. Let us make no comment here about when this
condition is sufficient. But we must point out that the object
$\ft{\delta_\Lambda}$ is not necessarily a measure but a so called {\em
tempered distribution} \cite{rudin1973fa}.

\subsection{A problem from the Scottish Book}\label{sec:scottish}

Nice reformulation but what could this buy for us?

Following the time-honored tradition of our
trade, let us first generalize: whatever we've said so far holds just as well
for an (almost) arbitrary function $f$ (nonnegative, integrable would do) in
place of $\chi_\Omega$. In other words, if one has
\beql{soft-tiling}
\sum_{\lambda \in \Lambda} f(x-\lambda) = \ell = {\rm Const.},
\eeq
for almost every $x \in \RR^d$ then it follows that
\beql{soft-ft-cond}
\supp \ft{\delta_\Lambda} \subseteq {\mathcal Z}(\ft{f})
\cup \Set{0} \mbox{\ \ where ${\mathcal Z}(\ft{f}) = \Set{\ft{f}=0}$}.
\eeq
Whenever \eqref{soft-tiling} holds we will say from now on that $f$ tiles space
with $\Lambda$ at level $\ell$.

Let us apply our new tool to Problem No 181 in the famous Scottish
Book\footnote{A collection of problems from the Scottish Caf\'e in Lw\'ow (then
in Poland, now Lviv in the Ukraine) where a number of important mathematicians
met for years exchanging problems. These problems, and whatever
solutions, were recorded in a notebook that formed the basis for
\cite{mauldin1981scottish}.}
\cite{mauldin1981scottish}, a problem posed
by H. Steinhaus and partly solved by him around 1939.
\begin{center}
 \resizebox{!}{7cm}{\includegraphics{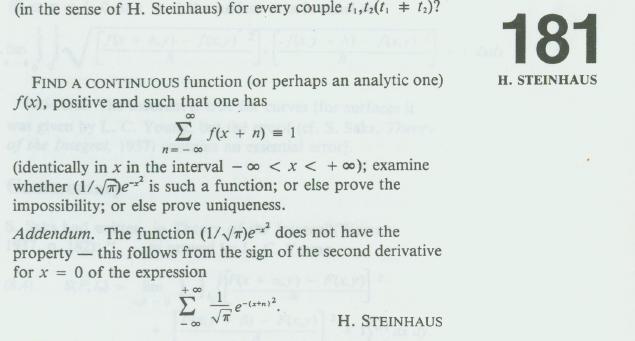}}
\end{center}
In our language the question essentially was if there is a strictly positive, continuous
and, if possible, analytic function $f$ that tiles $\RR$ with $\ZZ$.

Whoever has been exposed even a little to the FT knows that $e^{-x^2}$ has the
nice property that its FT transform is essentialy itself (times a constant and
possibly rescaled, depending on the normalizations used in the definition of the
FT) and therefore it has no zeros at all. It follows that \eqref{soft-ft-cond} has
no chance of holding as $\ft{\delta_\Lambda}$ has no place to be supported at
apart from 0, and it's not a constant. And this argument does not use the fact
that $\Lambda = \ZZ$ so $e^{-x^2}$ does not tile $\RR$ with {\em any} set of
translates.

An analytic solution to this question of Steinhaus can also be provided using
our technology. First of all, let us make the important remark that
\beql{psf}
\ft{\delta_{\ZZ^d}} = \delta_{\ZZ^d}.
\eeq
This is the deservedly famous Poisson Summation Formula
\cite[Ch.\ 9]{rudin1987real}
$$
\sum_n f(n) = \sum_n \ft{f}(n)
$$
in disguise.
In this case, when $\ft{\delta_\Lambda}$ is a measure, condition
\eqref{soft-ft-cond} is also sufficient for tiling
\cite{kolountzakis2004milano}, so it is sufficient to find a positive function
$f$ whose FT vanishes on $\ZZ\setminus\Set{0}$. Analyticity of $f$ will be
guaranteed if $\ft{f}$ has compact support.
It is then easy to check that one can take as $\ft{f}$ the sum of two triangle
functions with incommensurable base lengths,
\begin{center}
\resizebox{7cm}{2cm}{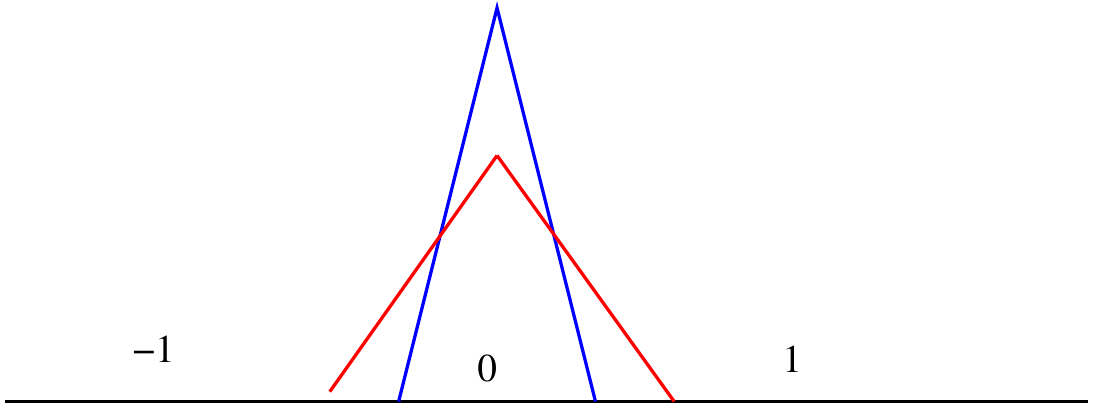}
\end{center}
as the FT of each such triangle is a nonnegative function (a {\em Fej\'er}
kernel) whose zeros are at integer multiples of the reciprocal of its half-base
length. Since these will never match for the two Fej\'er kernels their sum is always
positive. Since we also take the bases of the two triangles to be supported in
$(-1, 1)$ it follows that \eqref{soft-ft-cond} holds and we have tiling.

\subsection{Filling a box with two types of bricks}\label{sec:bricks}

Let us give another amusing application of the FT method to tiling problems. This will be a Fourier analytic proof \cite{kolountzakis2004box} of a result of Bower and Michael \cite{bower2004box}. Suppose that you have two types of rectangular bricks at your disposal, type $A$ of dimensions $a_1 \times a_2$ and type $B$ of dimensions $b_1 \times b_2$ (we're stating everything in dimension 2 but everything works in any dimension) and your task is to tile a rectangular box, say $Q=(-1/2, 1/2)^2$, with copies of bricks A and B.
The bricks may be translated but not rotated.

We will show that this is possible if and only if you can cut $Q$ along the $x$ or along the $y$ direction into two rectangles each of which can be tiled using bricks of one type only. A generic tiling of the type shown below left implies the existence of a {\em separated} tiling as shown below right.

\begin{center}
 \resizebox{10cm}{!}{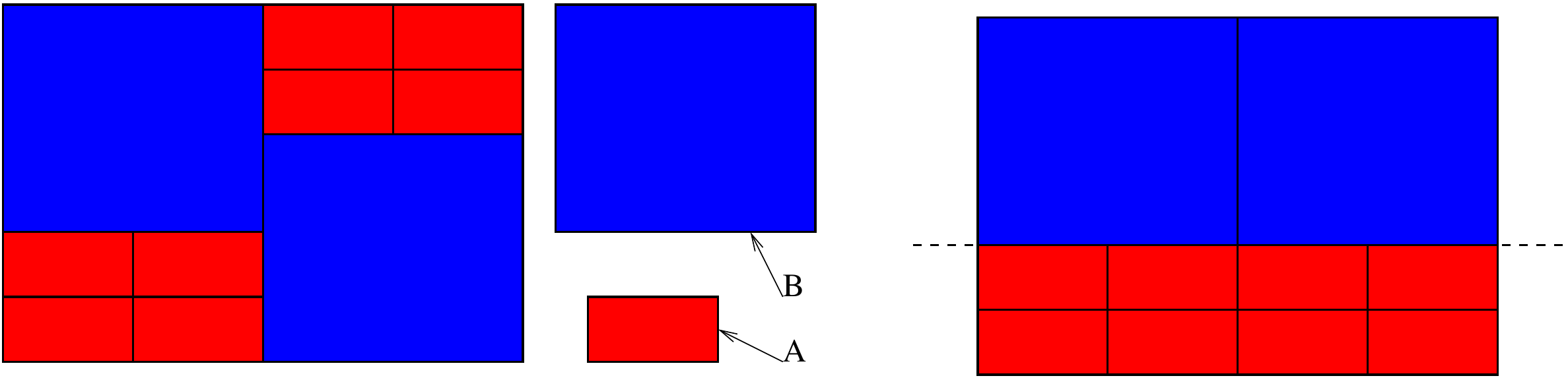}
\end{center}

\noindent
In other words, having two types of bricks at your disposal does not
demand much ingenuity on your side in order to exploit for tiling purposes.
If you cannot do it in an obvious way (i.e.\ without mixing the two types of bricks)
then it cannot be done at all.

A simple calculation shows that if $C = (-\frac{c_1}{2}, \frac{c_1}{2}) \times (-\frac{c_2}{2}, \frac{c_2}{2})$ is a centered $c_1 \times c_2$ box then
$$
\ft{\chi_C}(\xi, \eta) = \frac{\sin(\pi c_1 \xi)}{\xi} \cdot \frac{\sin(\pi c_2 \eta)}{\eta},
$$
and therefore $\ft{\chi_C}$ vanishes at those points, shown below as solid lines,

\def\labelone{$1/c_1$}
\def\labeltwo{$1/c_2$}
\begin{center}
 \resizebox{7cm}{!}{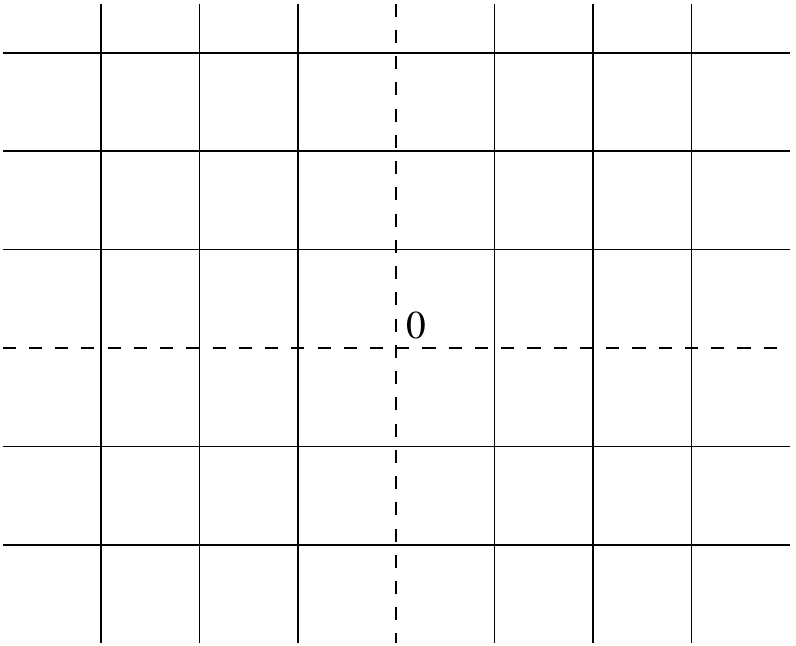}
\end{center}

\noindent
where the $\xi$ or the $\eta$ coordinate is a non-zero multiple of $1/c_1$ or $1/c_2$ respectively

Suppose now that we can tile $Q$ by translating copies of brick $A$ to the locations $T$ and copies of brick $B$ to the locations $S$. In other words
\beql{box-filled}
\chi_Q(x) = \sum_{t \in T} \chi_A(x-t) + \sum_{s \in S} \chi_B(x-s).
\eeq
which we rewrite as
$$
\chi_Q = \chi_A * \delta_T + \chi_B * \delta_S
$$
and take the FT of both sides to get
\beql{ft-box-filled}
\ft{\chi_Q}(\xi, \eta) = \phi_T(\xi, \eta)\ft{\chi_A}(\xi, \eta) + \phi_S(\xi, \eta) \ft{\chi_B}(\xi, \eta),
\eeq
where $\phi_T = \ft{\delta_T}$ and $\phi_S = \ft{\delta_S}$ are two {\em trigonometric polynomials}. Since $Q = (-1/2, 1/2)^2$ we have
\beql{zero-set-q}
{\mathcal Z}(\ft{\chi_Q}) = \Set{\ft{\chi_Q}=0} =
\Set{(\xi, \eta):\ \xi \in \ZZ\setminus\Set{0} \mbox{\ or\ } \eta \in \ZZ\setminus\Set{0}}.
\eeq
But, because of \eqref{ft-box-filled}, $\ft{\chi_Q}$ must vanish on the common zeros of $\ft{\chi_A}$ and $\ft{\chi_B}$, for instance at the points $(1/a_1, 1/b_2)$ and $(1/b_1, 1/a_2)$, which implies, because of \eqref{zero-set-q}
\beql{cases}
(1/a_1 \in \ZZ \mbox{\ or\ } 1/b_2 \in \ZZ) \mbox{\ \ and\ \ } (1/b_1 \in \ZZ \mbox{\ or\ } 1/a_2 \in \ZZ).
\eeq
If \eqref{cases} is true because $1/a_1, 1/a_2 \in \ZZ$ then brick $A$ alone can fill $Q$. Similarly if \eqref{cases} is satisfied with $1/b_1, 1/b_2 \in \ZZ$ then brick $B$ alone suffices.

What happens if $1/a_1, 1/b_1 \in \ZZ$? Since we have assumed a tiling of $Q$ in \eqref{box-filled} it follows, by traversing the $y$-axis, that there are nonnegative integers $k, l$ such that
$$
1 = k a_2 + l b_2.
$$
Cut now the $Q$ box parallel to the $x$-axis at height $ka_2$ from the bottom as show here:

\begin{center}
 \resizebox{3cm}{!}{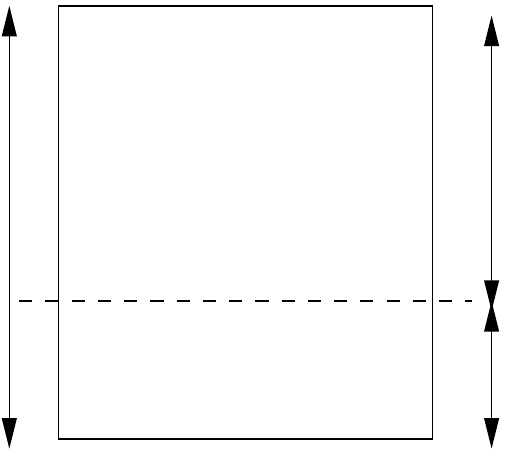}
\end{center}

\noindent
Now it is clear that brick $A$ can fill the lower box (since $1/a_1 \in \ZZ$) and brick $B$ can fill the upper box (since $1/b_1 \in \ZZ$). The remaining case $1/a_2, 1/b_2 \in \ZZ$ is treated similarly.

\section{Discrete tilings}\label{sec:discrete}

\subsection{Tilings of the integers and periodicity}\label{sec:tiling-the-integers}
Let us now focus on tiling the integers.
Let $A\subseteq \ZZ$ be a finite set and $\Lambda\in\ZZ$. We say that $A$ tiles
$\ZZ$ with $\Lambda$ at level $\ell$ if the copies $A+\lambda$, $\lambda \in
\Lambda,$ cover every integer exactly $\ell$ times. In other words
$$
\sum_{\lambda\in\Lambda} \chi_A(x-\lambda)=\ell, \mbox{\ \ for all $x \in \ZZ$}.
$$
We denote this situation as $A+\Lambda = \ell\ZZ$.

We say that a tiling is periodic with period $t\in\ZZ$ if $\Lambda+t = \Lambda$.
It is a basic fact proved by Newman \cite{newman1977tesselations} that all tilings of the integers at level 1 are periodic.

Indeed, suppose that $A=\Set{0=a_1<\ldots<a_k}$ is a finite set (we may freely
translate $A$ without changing the tiling or the periodicity property so we
assume it starts at $0$) and $A+\Lambda = \ZZ$ is a tiling at level 1. Fix any
$x\in\ZZ$ and write $W_x={x, x+1, \ldots, x+a_k-1}$ for the ``window'' of width
$a_k-1$ (one less than $A$) starting at $x$. We claim that the set $\Lambda
\cap W_x$ determines $\Lambda$. Let us show that it determines $\Lambda$ to the
right of $x+a_k-1$. It will determine $\Lambda$ to the left of $x$ by the same
argument.

It is enough to decide, looking at $\Lambda \cap W_x$ only, if $x+a_k \in
\Lambda$ or not. (We then repeat for $x+a_k+1$ and so on.) Observe that for any
$\lambda \in \Lambda \cap (-\infty, x)$ the set $A+\lambda$ is contained in
$(-\infty,x+a_k-1]$, so any such copy cannot be used to cover $x+a_k$. Clearly
no copy of the form $A+\lambda$ for $\lambda>x+a_k$ can be used for that
purpose too. We conclude that $x+a_k$ is covered by some copy $A+\lambda$ with
$\lambda \in W_x \cup \Set{x+a_k}$. Inspecting $\Lambda\cap W_x$ we can tell if
the relevant $\lambda$ is in $W_x$ or not. If it is in $W_x$ then $x+a_k \notin
\Lambda$ since this would lead to the copies $A+\lambda$ and $A+(x+a_k)$ to
overlap at $x+a_k$. If it is not in $W_x$ then necessarily $x+a_k \in \Lambda$,
and this concludes the proof of the claim.

How many different values can the set $\Lambda \cap W_x$ take? Clearly it can
take at most $2^{a_k}$ different values as there are two choices (in $\Lambda$
or not in $\Lambda$) for each $x \in W_x$. This means that there are two
different $x, y \in \Set{0,1,\ldots,2^{a_k}}$ for which $\Lambda \cap W_x$ is a
translate of $\Lambda \cap W_y$. It follows that $\Lambda + (y-x) = \Lambda$.
We have proved that every tiling has a period which is at most $2^D$ where $D$
is the diameter of the tile.

There is a similar result for tilings of the continuous line by translates of a
function
\cite{kolountzakis1996structure,lagarias1996tiling,kolountzakis2004milano}.
Combinatorial arguments do not seem to be enough here and the Fourier analytic
technology along with some deep results of Harmonic Analysis are being used for
the proof.

\subsection{Tilings of the finite cyclic group}\label{sec:cyclic}
The fact that a tiling of the integers is periodic allows us to view it as a
tiling on a smaller structure, a cyclic group. Indeed, assume that the tiling
$A+\Lambda=\ZZ$ has period $n$, that is $\Lambda+n=\Lambda$. Define the set
$$
\tilde\Lambda = \Lambda \bmod n \subseteq \Set{0,1,\ldots,n-1}
$$
by taking for each $\lambda \in \Lambda$ its residue $\bmod\ n$. It follows
from the $n$-periodicity of $\Lambda$ that $\Lambda = \tilde\Lambda + n\ZZ$
(this is again a tiling or a {\em direct sum}: every element of $\Lambda$ can
be written in a unique way as an element of $\tilde\Lambda$ plus an element of
$n\ZZ$).  Hence we have
$$
\ZZ = A + \Lambda = A + \tilde\Lambda + n\ZZ
$$
with all sums being direct. Taking quotients we obtain that the cyclic group
$\ZZ_n = \ZZ/(n\ZZ)$ can be written as a direct sum (tiling)
\beql{cyclic-group-tiling}
\ZZ_n = A + \tilde\Lambda.
\eeq
In this case we obviously have $n = \Abs{A}\cdot\Abs{\tilde\Lambda}$.

Let us stop here to make two side remarks about periodicity. The first is that
tilings of the cyclic group $\ZZ_n$ can also be periodic. Indeed, assume $\ZZ_n
= A + \tilde\Lambda$ and there exists $k\in \ZZ_n\setminus\Set{0}$ such that $A$ is
periodic with period $k$, i.e. $A+k=A$. Then we can reduce the set $A$ modulo $k$,
$\tilde A = A \bmod k \subseteq \Set{0,1,\ldots,k-1}$, and conclude as above
that $\ZZ_k = \tilde A + \tilde\Lambda$. Therefore, periodic tilings of $\ZZ_n$
can be regarded as tilings of a smaller group $\ZZ_k$ (it is trivial to see that
$k$ automatically divides $n$). It is thus natural to ask whether certain cyclic
groups $\ZZ_n$ admit {\em only} periodic tilings, i.e. whenever $\ZZ_n = A +
\tilde\Lambda$, then either $A$ or $\tilde\Lambda$ {\em must} be periodic. These
were called ``good'' groups by Haj\'os \cite{hajos1950factorization} (the
notion also makes sense in the more general setting of finite Abelian groups,
not only cyclic groups).  It turns out that some groups indeed have this
property, and Sands completed the classification of good groups in
\cite{sands1957factorisation,sands1962factorisation}. In particular, the good
groups that are cyclic are $\ZZ_n$ where $n$ divides one of $pqrs$, $p^2qr$,
$p^2q^2$ or $p^mq$, where $p, q, r, s$ are any distinct primes. The cyclic
group of the smallest order which is not good is $\ZZ_{72}$.

The other remark concerns the connection of tilings of cyclic groups to {\em
music composition}. As explained above, periodic tilings of $\ZZ_n$ are
mathematically speaking less interesting because they can be considered as
tilings of some smaller group $\ZZ_k$. It turns out that non-periodic tilings
of a cyclic group $\ZZ_n$ are also more interesting from an {\em aesthetic}
point of view, and they are called {\em Vuza-canons} in the musical community.
The interaction of mathematical theory and musical background has been
rather intense in recent years
\cite{andreatta2002tiling,amiot2004canons,andreatta2004group,fripertinger2004tiling,vuza1993supplementary}.
Finding  all non-periodic tilings of $\ZZ_n$ is therefore motivated by
contemporary music compositions. Of course, the problem makes sense only if
$\ZZ_n$ is not a good group (otherwise all tilings are periodic). Fripertinger
\cite{fri} achieved this task for $n = 72, 108$ while the authors
\cite{kolountzakis2009algorithms} gave an efficient algorithm to settle the
case $n = 144$ (the algorithm is likely to work for other values like $n = 120,
180, 200, 216$, beyond which the task simply seems hopeless).

Let us now return to the analysis of tilings of the cyclic group $\ZZ_n$. The Fourier condition \eqref{ft-cond} for tiling takes exactly the same form here and is much simpler to prove as no subtle analysis is required (no integrals, only finite sums are involved).
The Fourier transform of a function $f:\ZZ_n\to\CC$ is defined as the function $\ft{f}:\ZZ_n\to\CC$ given by
$$
\ft{f}(k) = \sum_{j=0}^{n-1} f(j) \zeta_n^{-kj},
$$
where $\zeta_n = e^{2\pi i/n}$ is a primitive $n$-th root of unity.

In the finite case the roles of the tile $A$ and the set of translations $\tilde\Lambda$ in \eqref{cyclic-group-tiling} can be interchanged, so let us adopt a more symmetric notation
$$
\ZZ_n = A + B,
$$
in which $A, B \subseteq \ZZ_n$ are just two subsets, necessarily satisfying $\Abs{A}\cdot\Abs{B} = n$, and such that every element of $\ZZ_n$ can be written uniquely as a sum of an element of $A$ and en element of $B$. The Fourier condition now takes the form
\beql{cyclic-ft-cond}
\ZZ_n = \Set{0} \cup {\mathcal Z}(\ft{\chi_A}) \cup {\mathcal Z}(\ft{\chi_B}).
\eeq
As we saw in the section \S\ref{sec:tiling-the-integers} every tiling of $\ZZ$
has a period which is at most $2^D$, where $D$ is the diameter of the tile. We
will see in this section how the easy combinatorial argument of
\S\ref{sec:tiling-the-integers} can be replaced with an argument that has the
Fourier condition \eqref{cyclic-ft-cond} as a starting point, uses some
well-known number-theoretic facts (properties of cyclotomic polynomials) and
gives much better results. We will describe the result in
\cite{kolountzakis2003long} and by I. Ruzsa in an appendix in
\cite{tijdeman2006periodicity} since it is simpler than the current best
results in \cite{biro2005divisibility}.

So, suppose that $A \subseteq \Set{0,1,\ldots,D}$ is a set of integers of
diameter $\le D$ and that $A$ tiles $\ZZ$ with period $M$. This implies that
$A$ (to be precise, $A$ reduced $\bmod\ M$) tiles the cyclic group $\ZZ_M$
$$
\ZZ_M = A+B.
$$
Assume also that $M$ is the least period. This implies that if $g \in \ZZ_M$ is
such that $B = B+g$ then $g=0$. (Otherwise the original tiling would be
periodic with period $g$.)

Now we use the Fourier transform without really mentioning it, using instead
polynomial terminology. Let
$$
A(x) = \sum_{a \in A} x^a,\ \ B(x) = \sum_{b \in B} x^b,
$$
be the polynomials defined by $A$ and $B$.
It is very easy to see that the tiling condition $\ZZ_M = A+B$ is the same as
$$
A(x)B(x) = 1 + x + x^2 + \cdots x^{M-1} \bmod\ x^M-1,
$$
or, in other words, that
\beql{divisibility}
x^M-1 \mbox{\ divides\ } A(x)B(x) - \frac{x^M-1}{x-1},
\eeq
and this implies that all $M$-th roots of unity except 1 are roots of the product $A(x)B(x)$.
(This last statement is the equivalent of \eqref{cyclic-ft-cond} in algebraic language.)

The $M$-th roots of unity $\zeta_M^j = e^{2\pi i j/M}, j=0,1,\ldots,M-1$, are
grouped into {\em cyclotomic classes}. Two roots $\zeta_M^j, \zeta_M^k$ belong
to the same cyclotomic class if and only if the greatest common divisors
$(j,M)$ and $(k,M)$ are the same. If a root of unity in a certain cyclotomic
class is a root of a polynomial with integer or rational coefficients then so
are all the other roots in the same cyclotomic class. If $d$ is a positive
integer then the $d$-th cyclotomic polynomial is the monic polynomial which has
as roots all the {\em primitive} $d$-th roots of unity
$$
\Phi_d(x) = \prod_{1\le j<d, (j,d)=1} (x-e^{2\pi i j/d}).
$$
We have $\deg \Phi_d(x) = \phi(d)$ (the Euler function) and we will use the estimate
\cite[p.\ 267]{hardy1979numbers}
\beql{estimate}
C\frac{d}{\log\log d} \le \phi(d) \le d,
\eeq
which holds for some constant $C$ and all large enough $d$.  It turns out that
each $\Phi_d$ is an inrreducible polynomial with integer coefficients and one
can write
$$
x^M-1 = \prod_{d\mbox{\tiny \ divides\ }M} \Phi_d(x).
$$
Coupled with \eqref{divisibility} this implies that if $d>1$ is a divisor of
$M$ then $\Phi_d(x)$ divides the product $A(x)B(x)$.

Let now $\Phi_{s_1}(x),\ldots,\Phi_{s_k}(x)$ be all cyclotomic polynomials
$\Phi_s(x)$ with $s>1$ that divide $A(x)$, written {\em once} each
and numbered so that $1<s_1<s_2<\cdots<s_k$.
Since $\deg\Phi_s = \phi(s)$ it follows that
\beql{degree-inequality}
\phi(s_1)+\cdots+\phi(s_k) \le \deg{A(x)} \le D.
\eeq
We have
\begin{alignat*}{3}
\sum_{i=1}^k s_i & \le C \sum_{i=1}^k \phi(s_i)\log\log s_i & \mbox{(from \eqref{estimate})} \ \\
 & \le C \sum_{i=1}^k \phi(s_i)\log\log D & \qquad\\
 & \le C D \log\log D & \mbox{(from \eqref{degree-inequality})}.
\end{alignat*}
From this inequality and using the fact that the $s_j$ are {\em different}
integers it follows that
\beql{k-bound}
k \le C \sqrt{D \log\log D},
\eeq
as obviously the worst case is if the $s_j$ are $k$ consecutive integers and in
that case their sum is $\gtrsim k^2$.

Define now $t = \prod_{i=1}^k s_i$,
so that all cyclotomic polynomials that divide $A(x)$ are also
divisors of $x^t-1$.
It follows that
$$
x^M-1 \mbox{\ divides\ } (x^t-1)B(x)
$$
and this means precisely that $B$ has period $t$, hence $t\ge M$ as we assumed $M$ to be minimal.
Using \eqref{k-bound} and the bound
$s_j = O(D^2)$, for instance, we get
the following bound for the least period $M$:
\beql{t-bound}
M \le t \le \exp(C \sqrt{D} \log D \sqrt{\log\log D}),
\eeq
for an appropriate constant $C$. This is a much better bound, in terms of its dependence on $D$, that the bound $M \le 2^D$ that we obtained with combinatorics alone.

\subsection{Periodicity in two dimensions and computability}\label{sec:computability}
In contrast to the one-dimensional case, where every tiling by a finite subset
of $\ZZ$ is periodic, it is very easy to see that in $\ZZ^2$ there are tilings
which are not periodic. One has to be a little careful with periodicity in
dimension 2 and higher though. We shall call a set $A \subseteq
\ZZ^d$
periodic if there exists a full lattice of periods, i.e.\ if there exist $d$
linearly independent vectors $u_1,\ldots,u_d$ which are all periods of the set
$$
A + u_1 = A + u_2 = \cdots = A + u_d = A.
$$
If that happens then every integer linear combination of the $u_j$ is also a
period or, in other words, the set $A$ has a full-dimensional lattice of periods
$$
\Lambda = {\rm span}_{\ZZ}\, \Set{u_1, u_2, \ldots, u_d}.
$$
Thus in dimension $2$ a set may have a period but not two linearly independent
periods, in which case we do not call it periodic. One example is the set
$\ZZ\times\Set{0} \subseteq \ZZ^2$.
A tiling of $\ZZ^2$ by a finite subset which is not periodic is very easy to
construct.
\begin{center}
 \resizebox{3cm}{!}{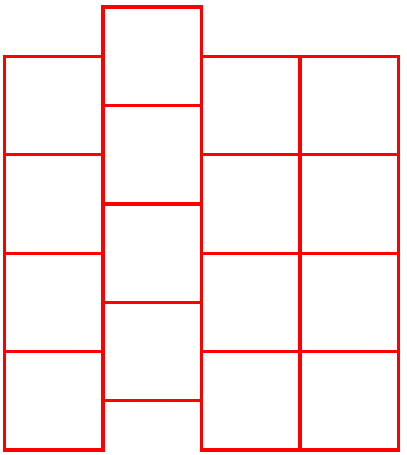}
\end{center}
On takes $A$ to be a square, for instance the set $\Set{(0,0), (1,0), (1,0),
(1,1)}$ and perturbs the usual tiling of it by shifting one column only up by
one unit. This destroys the period $(2,0)$ of the tiling along the $x$-axis,
but leaves the period $(0,2)$ along the $y$-axis intact.

It is not entirely obvious how to construct a tiling of $\ZZ^2$ which has no
period at all, but it can be done and let us briefly describe how. The key idea
is to start again with the usual square tiling and do something to it so as to
destroy all periods, not just the periods along one axis. The way to do this is
to simultanesouly shift a horizontal and a vertical column. This is not of
course possible with the square tile we used before, as simple experimentation
will convince you. But it can be done if we ``interleave'' four tilings of the
type shown above. Take our tile to be the set
$$
A = \Set{(0,0), \ (2,0), \ (0,2), \ (2,2)}
$$
which is what we had before only dilated by a factor of two. Now tile
$(2\ZZ)^2$ using $A$ in the usual way, that is by translating $A$ to the
locations $(4\ZZ)^2$. Using this tiling one can get a tiling of $\ZZ^2$ by
tiling the four sets (cosets of the subgroup $(2\ZZ)^2$ in $\ZZ^2$)
$$
(2\ZZ)^2, \ (2\ZZ)^2+(1,0), \ (2\ZZ)^2+(0,1), \ (2\ZZ)^2+(1,1)
$$
in exactly the same way. But the smart way to do it is to use the nice tiling
for the first two of the above four cosets, then use, for tiling
$(2\ZZ)^2+(0,1)$, a tiling like that shown in the figure above (which destroys
the horizontal periods) and finally use for tiling $(2\ZZ)^2+(1,1)$ a similar
tiling which destroys the vertical periods.
The tiling that we constructed in this manner has no period vectors at all.

Having established that in two dimensions there are tilings with no periods at
all, let us now remark that it is very different to ask for {\em tiles} which
are aperiodic, i.e.\ for tiles which can tile but only in a way that admits no
periods. In fact, the answer to this question is not known if we insist that
only translations are allowed.

\begin{conjecture}[Lagarias and Wang \cite{lagarias1996tiling}]
If a finite subset $A \subseteq \ZZ^2$ can tile by translation then it can also
tile in a periodic way.
\end{conjecture}

Let us point out that if more freedom than translation is allowed then tiles
(or sets of tiles) are known which can only tile aperiodically, so the above
conjecture, if true, would mean that restricting the allowed motions to
translations makes a big difference in this respect.

Let us also remark that it can be proved, at least in dimension 2
\cite{rao2006tiling}, that if a set admits a
tiling with just one period then it also admits a (fully) periodic tiling.

The property of periodicity is of great importance to the issue of
computability, namely to whether a computer can decide if a finite set $A$
admits tiling by translation. In still more words, we are interested to know if
there is a computer algorithm (a Turing machine for purists) which, given a
finite subset $A$ of $\ZZ^2$, will decide in finite time if there exists a
tiling complement of $A$ or not, i.e.\ if there exists $B \subseteq \ZZ^2$ such
that $A+B=\ZZ^2$ is a tiling.

We do not require our algorithm to:
\begin{itemize}
 \item run fast (only to finish at some point)
 \item find such a tiling complement $B$ if it exists.
\end{itemize}
Especially for the second point notice that it would not make sense to want to
know $B$ as such a set is an infinite object and there is no a priori reason
for it to be describable in some finite manner.

It is easy to give an algorithm that would answer NO if $A$ is not a tile but
would run forever if $A$ is a tile. Here is how this can be done. For each
$n=1,2,\ldots$, decide by trial and error if there is a finite set $B$ such that
$A+B$ covers the (discrete) square $Q_n = [-n,n]^2 \cap \ZZ^2$ in a
non-overlapping way. There is a finite number of such sets $B$ to try as it does
not matter what $B$ is far away from $Q_n$. If such a $B$ exists move on to the
next $n$. If not then declare $A$ a non-tile and stop. The correctness of this
algorithm (i.e., the fact that the algorithm will stop for any non-tile $A$)
follows from a fairly simple diagonal argument: if a finite
$A$ can cover in a nonoverlapping
way every $Q_n$ then
it can also tile $\ZZ^2$. We invite the reader to think over why this is so.

The hard thing of course is to have an algorithm that always halts and say YES
if $A$ is a tile and NO if it is not a tile. Let us point out here that in
various other tiling situations such an algorithm does not exist. For instance
\cite{golomb1970tiling} there is no algorithm which, given
a finite collection of subsets of $\ZZ^2$, decides if these can tile the
plane by translation. This remains true even if the number of subsets remains
bounded (but large).

The connection with periodicity \cite{robinson1971undecidability} is that
the periodic tiling conjecture above implies decidability. Indeed, let us
assume that any set that tiles can also tile periodically. (We emphasize here
that we are making no assumption about the size of the periods or their
dependence on the tile.) This means that it is equivalent to decide if a given
set $A$ admits periodic tilings.

Let us assume that $A \subseteq [0, D]^2$. For $n=1,2,\ldots$, we can clearly
find all subsets $B$ of $[-n-2D, n+2D]^2$ such that $A+B$ is a nonoverlapping
covering of $Q_n = [-n, n]^2$. We can do this, slowly but surely, by just
trying all eligible subsets. If we find none then clearly our set is not a
tile, we say NO, and we stop. If we do find some we look among them to find a
set $B$ which ``can be extended periodically''. How does such a set look like?

If an infinite set $\tilde{B} \subseteq \ZZ^2$ is periodic with linearly
independent period vectors $u_1, u_2 \in \ZZ^2$ then we can always find two
vectors
$$
\tilde{u}_1 = (a, 0),\ \tilde{u}_2 = (0,b),\ \ a, b \in \Set{1,2,3,\ldots},
$$
which are also periods (in other words, every lattice in $\ZZ^2$ contains a
sublattice generated by two vectors on the $x$- and $y$-axes).
The mapping $(x,y) \to (x \bmod a, y \bmod b)$ maps $\tilde{B}$ to a set
$$
B' \subseteq \Set{0,1,\ldots,a-1}\times\Set{0,1,\ldots,b-1}
$$
and the periodicity of $\tilde{B}$ means exactly that
$$
\tilde{B} = (a\ZZ)\times(b\ZZ) + B'
$$
is a direct sum, i.e.\ every element of $\tilde{B}$ can be written uniquely in
the form $(am, bn) + b'$, where $m, n \in \ZZ, b' \in B'$.

Suppose now that our algorithm looks at a finite set $B$ which is a finite part
of a periodic set $\tilde{B}$. This finite part of $\tilde{B}$ arises if we
keep those elements of $\tilde{B}$ which are used in a non-overlapping covering
of $Q_n$ (shown as a dashed green rectangle). If $n$ is large enough then the
set $B$ will look something like this:

\begin{center}
 \resizebox{7cm}{!}{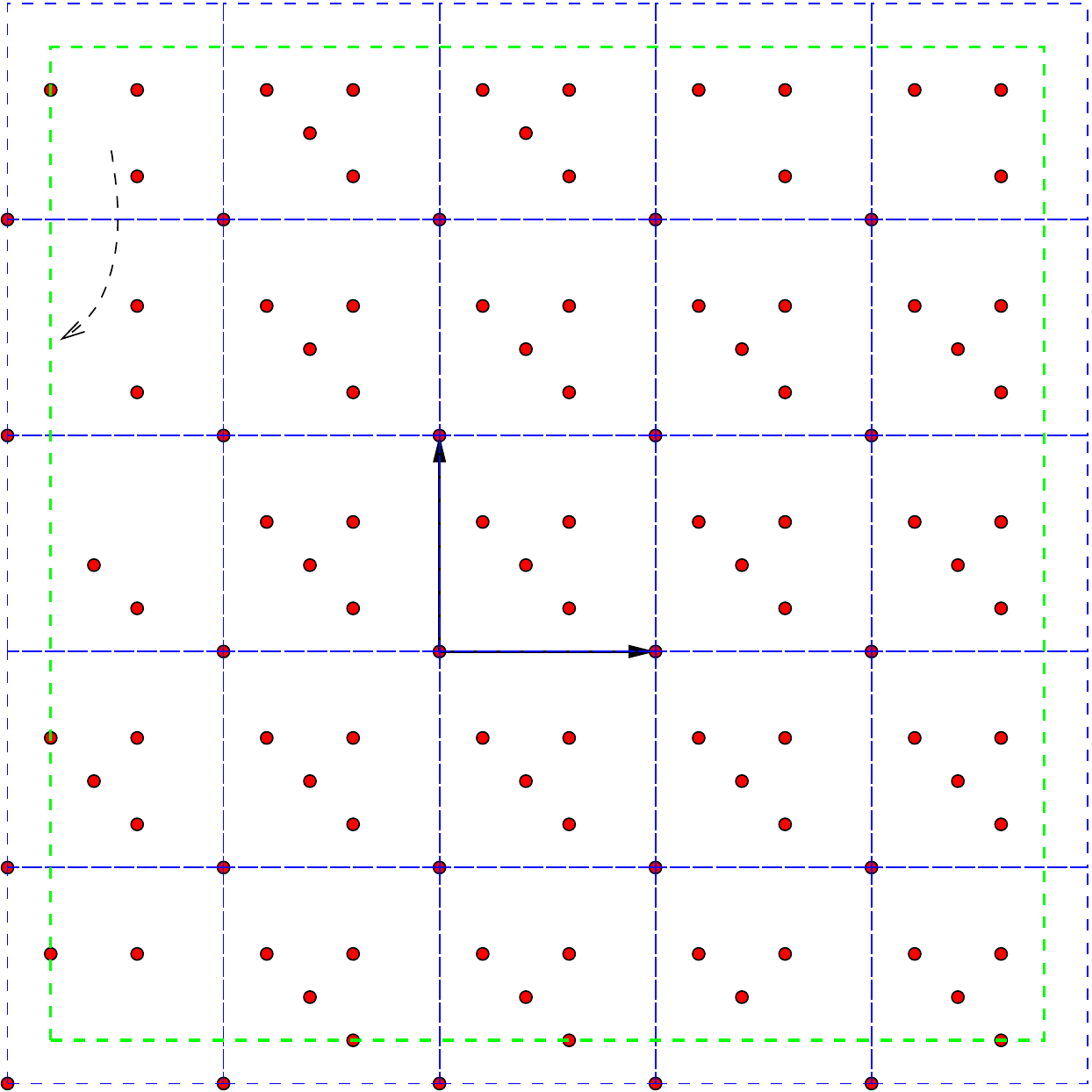}
\end{center}

That is, the set $B$ will contain several ``full periods'' (the dashed blue
$a\times b$ rectangles) $\tilde{B}$ plus a few incomplete periods near the
border of $Q_n$. Suppose now that we have a $k \times k$
block of full periods in the set $B$ and that $k$ is large enough that both
$ka$ and $kb$ are much larger than the diameter of $A$, say larger than $10
D$. Then we can just extend this block periodically (with periods $(a,0)$
and $(0,b)$) and we will get a periodic tiling complement of $A$.
The reason is that (a) no overlaps arise in this manner or they would have
shown up in the $k \times k$ block already and that (b) we have covering of
everything as clearly the $k \times k$ block suffices to cover the square of
side $5D$ with the same center.

Thus our algorithm will, sooner or later, either discover that $A$ cannot tile
(since it will not be able to tile a finite region) or it will discover that $A$ can tile
a large region with a tiling complement that looks like the one in the figure above,
and which can therefore be extended periodically to a full tiling of $\ZZ^2$.
This means that our algorithm will always stop and give us a correct YES or a correct NO.
We have not given any bound on the running time though.

\section{Tiles, spectral sets and complex Hadamard matrices}

In this last section we introduce the notion of spectral sets and  describe how
it is related to tilings and complex Hadamard matrices. For the sake of
simplicity, we will restrict our attention to finite groups, in particular to
cyclic groups $\ZZ_n$ and its powers $\ZZ_n^d$.

Although we have seen in the sections above that many fundamental  problems
about translational tilings remain open, the reader will agree with us that the
notion of tiling is very intuitive and easy to grasp. To the contrary, the
definition of spectral sets is somewhat more abstract. In order to make it more
down-to-earth, we will not use the standard definition here, but first introduce
complex Hadamard matrices and then use them to define spectral sets.

What are Hadamard matrices? Classically, they are square matrices consisting of
elements $\pm 1$ only, such that the rows (and thus also the columns) are
orthogonal to each other. Complex Hadamard matrices are a natural generalization
of this concept. A $k\times k$ matrix $H$ is a complex Hadamard matrix if its
entries are complex numbers of modulus 1, and the rows (and thus the columns)
are orthogonal. Recall here that orthogonality is understood with respect to the
{\em complex} scalar product, i.e. you need to conjugate in one of the
components, the scalar product of $(z_1, \dots z_k)$ and $(u_1, \dots u_k)$
being $\sum_{j=1}^k \overline{z_j}u_j$.

It is well-known that a $k\times k$ real Hadamard matrix can only exist if $k$
is divisible by 4. On the contrary, complex Hadamard matrices exist in all
dimension. Indeed the {\it Fourier matrix} $F_k$ defined as
$$F_k:=\left [
\begin{array}{cccccc}
1&1&\cdot&\cdot&\cdot&1\\
1&\zeta_k&\zeta_k^2&\cdot&\cdot&\zeta_k^{k-1}\\
\cdot&\cdot&\cdot&\cdot&\cdot&\cdot\\
\cdot&\cdot&\cdot&\zeta_k^{jm}&\cdot&\cdot\\
\cdot&\cdot&\cdot&\cdot&\cdot&\cdot\\
1&\zeta_k^{k-1}&\cdot&\cdot&\cdot&\zeta_k
\end{array}
\right ]$$

where $\zeta_k = e^{2\pi i/k}$, is a complex Hadamard matrix for every $k$.

Each element $h_{j,m}$ of a $k\times k$ complex Hadamard matrix $H$ is of
the form $e^{2i\pi \rho_{j,m}}$ where $\rho_{j,m}\in [0,1)$. We will call the
real $k\times k$ matrix $R$ formed by the angles $\rho_{j,m}$ the (entry-wise) {\it
logarithm} of $H$. In notation, $R=\frac{1}{2i\pi} {\rm LOG} (H)$ or $H={\rm EXP} (2i\pi R)$ (note that both these operations are meant entry-wise). For a real matrix
$R$ such that ${\rm EXP} (2i\pi R)$ is complex Hadamard, we will use the terminology that
$R$ is a {\it log-Hadamard} matrix.

Complex Hadamard matrices play an important role  in quantum information theory,
in particular in the construction of teleportation and dense coding schemes
\cite{werner2001teleportation}, among other things. An online catalogue of all
known families of complex Hadamard matrices is available at
\cite{website2010hadamards}.

Let us now turn to the definition of spectral sets. We identify elements of
$\ZZ_n^d$ with column vectors, each coordinate being in the range $\{0,
1, \dots , n-1\}$. Thus, a set $S\subset \ZZ_n^d$ with $r$ elements can be
identified with an $d\times r$ matrix, the columns of which are the elements of
$S$ (the order of elements does not matter). We will abuse notation and denote
this matrix also by $S$. We also identify $\ft{\ZZ_n^d}$ with row
vectors whose coordinates are in $\{0, 1, \dots , n-1\}$. Accordingly, we will
use the notation $\ft{\ZZ_n^d}=\ZZ_n^{d \top}$ (the $\top$ meaning
transposition). Two sets $S\subset\ZZ_n^d$ and $Q\subset \ZZ_n^{d \top}$, both
of them having $r$ elements, are called a {\it spectral pair} if the matrix
$\frac{1}{n}Q S$ is an $r\times r$ log-Hadamard matrix. In this case $Q$ is
called the {\it spectrum} of $S$, and $S$ is called a {\it spectral set}. We
remark that $S^\top Q^\top =(QS)^\top$ (where $S^\top$ and $Q^\top$ denote the
transposed matrices) is then automatically also a log-Hadamard matrix, so we
obtain that $S$ is a spectrum of $Q$, which justifies the symmetry of the
terminology "spectral pair".

For example, the whole group $\ZZ_n=\{0, 1, \dots , n-1\}$ is spectral, and its
spectrum is the whole dual group $\ZZ_n^{\top}$, giving rise to the logarithm of
the Fourier matrix, $\frac{1}{2i\pi} {\rm LOG} (F_n) = \frac{1}{n}\ZZ_n^\top \ZZ_n$. As a less trivial example, take the set 
$S=\{0,1,2\}$ in $\ZZ_6$.  Then $Q=\{0,2,4\}^{\top}$ is a spectrum of $S$, and $\frac{1}{6}QS=\frac{1}{2i\pi}{\rm LOG}(F_3)$, 
the logarithm of the $3\times 3$ Fourier matrix.

We saw that the tiling $\ZZ_n = A + B$ is characterized by $\Abs{A}\cdot\Abs{B}
= n$ and the Fourier condition $\ZZ_n = \Set{0} \cup {\mathcal Z}(\ft{\chi_A})
\cup {\mathcal Z}(\ft{\chi_B})$ in equation \eqref{cyclic-ft-cond}. After some
inspection one sees that spectral sets also admit a characterization by Fourier
analysis. Indeed, using the definition above we obtain that $S\subset \ZZ_n^d$,
$|S|=r$, is spectral if and only if there exists $Q\subset \ZZ_n^{d \top}$,
such that $|Q|=r$ and
\beql{spectral}
q_j-q_k\in {\mathcal Z}(\ft{\chi_S})
\eeq
for each $q_j\ne q_k\in Q.$ (This is an equivalent way of saying that the $j$th
and $k$th rows of the matrix ${\rm EXP} (\frac{2i\pi}{n}QS)$ are orthogonal.)

What are the connections of tiles and spectral sets? Notice that by the
orthogonality conditions the matrix $U=\frac{1}{\sqrt{r}}{\rm EXP} (\frac{2i\pi}{n}Q S)$
is unitary. This implies the following important fact (for a formal proof, see
\cite{kolountzakis2006tiles}):
\beql{power-spectrum-tiling}
\sum_{q \in Q} \Abs{\ft{\chi_S}}^2(x-q) = \Abs{S}^2=r^2,
\eeq
for all $x\in \ZZ_n^{d \top}$. In other words, the function
$\Abs{\ft{\chi_S}}^2$ tiles $\ft{\ZZ_n^d}$ at level $\Abs{S}^2$ when translated
at the locations $Q$. One can also say that the set $Q$ tiles the group
$\ft{\ZZ_n^d}$ with the {\em weighted translates} defined by the
nonnegative function
$$
\frac{1}{\Abs{S}^2}\Abs{\ft{\chi_S}}^2.
$$
The question thus
arises naturally: does $Q$ tile $\ft{\ZZ_n^d}$ in the ordinary sense, i.e. does
there exist a set $P\subset \ft{\ZZ_n^d}$ such that $Q+P=\ft{\ZZ_n^d}$? A famous
conjecture of Fuglede \cite{fuglede1974operators} concerns exactly this
scenario: it states that a set $S$ is spectral if and only if it tiles the group
(the conjecture was originally stated in the Euclidean space $\RR^d$ but it
makes sense in any abelian group).

Several positive partial results were proven concerning special cases of
Fuglede's conjecture (see e.g.
\cite{iosevich2001convexbodies,iosevich2003fuglede,kolountzakis2000nonsymmetric,
laba2001twointervals}), before T. Tao \cite{tao2004fuglede} showed an example of
a spectral set in $\RR^5$ which did not tile the space. Tao's example was based
on considerations in the finite group $\ZZ_3^5$. It also implies the existence
of such sets in any dimension $d\ge 5$. Subsequently, counterexamples in lower
dimensions (4 and 3, respectively) were found by the authors
\cite{matolcsi2005fuglede4dim,kolountzakis2006hadamard}. All these examples are
based on considerations in finite groups $\ZZ_n^d$ and, ultimately, on the
existence of certain complex Hadamard matrices.

The other direction of Fuglede's conjecture could not be settled by Tao's
arguments. Finally a non-spectral tile in $\RR^5$ was exhibited by the authors
\cite{kolountzakis2006tiles} by a tricky duality argument, also based on
considerations in finite groups. Later, counterexamples in  dimensions 4 and 3,
respectively, were found in \cite{farkas2006tiles,farkas2006onfuglede}. As of
today, Fuglede's conjecture is still open in both directions in dimensions 1 and
2.

As a final remark let us mention another interesting connection of tilings and complex Hadamard matrices. As discussed above, Fuglede's conjecture is not true in general, but it is true in many particular cases. What could this buy for us? We saw from the definition that spectral sets are directly related to complex Hadamard matrices, and the latter are very useful in another branch of mathematics: quantum information theory. Could we use the connection of tilings to spectral sets (which exists in many special cases) to construct new families of complex Hadamard matrices? It turns out that the answer is positive. It turns out \cite{matolcsi2007constructions} that "trivial" tiling constructions unfortunately lead to well-known families of complex Hadamard matrices (so-called Dita-families). However, a non-standard tiling construction of Szab\'o \cite{szabo1985atype} was used in \cite{matolcsi2007constructions} to produce previously unknown families of complex Hadamard matrices in dimensions 8, 12 and 16.

\bibliographystyle{abbrv}
\bibliography{ttsurvey}

\end{document}

%% file: one-square.pdf_t
\begin{picture}(0,0)%
\includegraphics{one-square.pdf}%
\end{picture}%
\setlength{\unitlength}{4144sp}%
\begingroup\makeatletter\ifx\SetFigFont\undefined%
\gdef\SetFigFont#1#2#3#4#5{%
  \reset@font\fontsize{#1}{#2pt}%
  \fontfamily{#3}\fontseries{#4}\fontshape{#5}%
  \selectfont}%
\fi\endgroup%
\begin{picture}(1374,1374)(2014,-2998)
\end{picture}%

%% file: cyclic.pdf_t
\begin{picture}(0,0)%
\includegraphics{cyclic.pdf}%
\end{picture}%
\setlength{\unitlength}{4144sp}%
\begingroup\makeatletter\ifx\SetFigFontNFSS\undefined%
\gdef\SetFigFontNFSS#1#2#3#4#5{%
  \reset@font\fontsize{#1}{#2pt}%
  \fontfamily{#3}\fontseries{#4}\fontshape{#5}%
  \selectfont}%
\fi\endgroup%
\begin{picture}(2581,106)(2873,-4839)
\end{picture}%

%% file: notched.pdf_t
\begin{picture}(0,0)%
\includegraphics{notched.pdf}%
\end{picture}%
\setlength{\unitlength}{3947sp}%
\begingroup\makeatletter\ifx\SetFigFont\undefined%
\gdef\SetFigFont#1#2#3#4#5{%
  \reset@font\fontsize{#1}{#2pt}%
  \fontfamily{#3}\fontseries{#4}\fontshape{#5}%
  \selectfont}%
\fi\endgroup%
\begin{picture}(3644,3344)(4179,-7583)
\end{picture}%

%% file: extended.pdf_t
\begin{picture}(0,0)%
\includegraphics{extended.pdf}%
\end{picture}%
\setlength{\unitlength}{3947sp}%
\begingroup\makeatletter\ifx\SetFigFont\undefined%
\gdef\SetFigFont#1#2#3#4#5{%
  \reset@font\fontsize{#1}{#2pt}%
  \fontfamily{#3}\fontseries{#4}\fontshape{#5}%
  \selectfont}%
\fi\endgroup%
\begin{picture}(3344,3644)(5079,-3083)
\end{picture}%

%% file: hexagon.pdf_t
\begin{picture}(0,0)%
\includegraphics{hexagon.pdf}%
\end{picture}%
\setlength{\unitlength}{4144sp}%
\begingroup\makeatletter\ifx\SetFigFont\undefined%
\gdef\SetFigFont#1#2#3#4#5{%
  \reset@font\fontsize{#1}{#2pt}%
  \fontfamily{#3}\fontseries{#4}\fontshape{#5}%
  \selectfont}%
\fi\endgroup%
\begin{picture}(5427,3681)(1336,-6880)
\put(1351,-5011){\makebox(0,0)[lb]{\smash{{\SetFigFont{12}{14.4}{\rmdefault}{\mddefault}{\updefault}{\color[rgb]{0,0,0}$\Omega$}%
}}}}
\put(3871,-6811){\makebox(0,0)[lb]{\smash{{\SetFigFont{12}{14.4}{\rmdefault}{\mddefault}{\updefault}{\color[rgb]{0,0,0}\pointslambda}%
}}}}
\end{picture}%

%% file: delta-lambda.pdf_t
\begin{picture}(0,0)%
\includegraphics{delta-lambda.pdf}%
\end{picture}%
\setlength{\unitlength}{4144sp}%
\begingroup\makeatletter\ifx\SetFigFont\undefined%
\gdef\SetFigFont#1#2#3#4#5{%
  \reset@font\fontsize{#1}{#2pt}%
  \fontfamily{#3}\fontseries{#4}\fontshape{#5}%
  \selectfont}%
\fi\endgroup%
\begin{picture}(7449,2734)(1564,-4798)
\end{picture}%

%% file: triangles.pdf_t
\begin{picture}(0,0)%
\includegraphics{triangles.pdf}%
\end{picture}%
\setlength{\unitlength}{4144sp}%
\begingroup\makeatletter\ifx\SetFigFont\undefined%
\gdef\SetFigFont#1#2#3#4#5{%
  \reset@font\fontsize{#1}{#2pt}%
  \fontfamily{#3}\fontseries{#4}\fontshape{#5}%
  \selectfont}%
\fi\endgroup%
\begin{picture}(4994,1844)(2229,-4583)
\end{picture}%

%% file: box-filled.pdf_t
\begin{picture}(0,0)%
\includegraphics{box-filled.pdf}%
\end{picture}%
\setlength{\unitlength}{4144sp}%
\begingroup\makeatletter\ifx\SetFigFont\undefined%
\gdef\SetFigFont#1#2#3#4#5{%
  \reset@font\fontsize{#1}{#2pt}%
  \fontfamily{#3}\fontseries{#4}\fontshape{#5}%
  \selectfont}%
\fi\endgroup%
\begin{picture}(10844,2609)(429,-2423)
\end{picture}%

%% file: zeros.pdf_t
\begin{picture}(0,0)%
\includegraphics{zeros.pdf}%
\end{picture}%
\setlength{\unitlength}{4144sp}%
\begingroup\makeatletter\ifx\SetFigFont\undefined%
\gdef\SetFigFont#1#2#3#4#5{%
  \reset@font\fontsize{#1}{#2pt}%
  \fontfamily{#3}\fontseries{#4}\fontshape{#5}%
  \selectfont}%
\fi\endgroup%
\begin{picture}(3624,2949)(1339,-3448)
\put(4681,-1951){\makebox(0,0)[lb]{\smash{{\SetFigFont{12}{14.4}{\familydefault}{\mddefault}{\updefault}{\color[rgb]{0,0,0}$\xi$}%
}}}}
\put(3241,-646){\makebox(0,0)[lb]{\smash{{\SetFigFont{12}{14.4}{\familydefault}{\mddefault}{\updefault}{\color[rgb]{0,0,0}$\eta$}%
}}}}
\put(3646,-2311){\makebox(0,0)[lb]{\smash{{\SetFigFont{12}{14.4}{\familydefault}{\mddefault}{\updefault}{\color[rgb]{0,0,0}\labelone}%
}}}}
\put(3196,-1546){\makebox(0,0)[lb]{\smash{{\SetFigFont{12}{14.4}{\familydefault}{\mddefault}{\updefault}{\color[rgb]{0,0,0}\labeltwo}%
}}}}
\end{picture}%

%% file: cut.pdf_t
\begin{picture}(0,0)%
\includegraphics{cut.pdf}%
\end{picture}%
\setlength{\unitlength}{4144sp}%
\begingroup\makeatletter\ifx\SetFigFont\undefined%
\gdef\SetFigFont#1#2#3#4#5{%
  \reset@font\fontsize{#1}{#2pt}%
  \fontfamily{#3}\fontseries{#4}\fontshape{#5}%
  \selectfont}%
\fi\endgroup%
\begin{picture}(2352,2049)(1894,-7408)
\put(2026,-6406){\makebox(0,0)[lb]{\smash{{\SetFigFont{12}{14.4}{\familydefault}{\mddefault}{\updefault}{\color[rgb]{0,0,0}1}%
}}}}
\put(4231,-6001){\makebox(0,0)[lb]{\smash{{\SetFigFont{12}{14.4}{\familydefault}{\mddefault}{\updefault}{\color[rgb]{0,0,0}$l b_2$}%
}}}}
\put(4231,-7081){\makebox(0,0)[lb]{\smash{{\SetFigFont{12}{14.4}{\familydefault}{\mddefault}{\updefault}{\color[rgb]{0,0,0}$k a_2$}%
}}}}
\end{picture}%

%% file: nonperiodic.pdf_t
\begin{picture}(0,0)%
\includegraphics{nonperiodic.pdf}%
\end{picture}%
\setlength{\unitlength}{4144sp}%
\begingroup\makeatletter\ifx\SetFigFont\undefined%
\gdef\SetFigFont#1#2#3#4#5{%
  \reset@font\fontsize{#1}{#2pt}%
  \fontfamily{#3}\fontseries{#4}\fontshape{#5}%
  \selectfont}%
\fi\endgroup%
\begin{picture}(1844,2069)(2004,-3233)
\end{picture}%

%% file: finite-periodic.pdf_t
\begin{picture}(0,0)%
\includegraphics{finite-periodic.pdf}%
\end{picture}%
\setlength{\unitlength}{4144sp}%
\begingroup\makeatletter\ifx\SetFigFont\undefined%
\gdef\SetFigFont#1#2#3#4#5{%
  \reset@font\fontsize{#1}{#2pt}%
  \fontfamily{#3}\fontseries{#4}\fontshape{#5}%
  \selectfont}%
\fi\endgroup%
\begin{picture}(5677,5676)(2661,-5050)
\put(4771,-2896){\makebox(0,0)[lb]{\smash{{\SetFigFont{12}{14.4}{\rmdefault}{\mddefault}{\updefault}{\color[rgb]{0,0,0}$0$}%
}}}}
\put(5986,-2671){\makebox(0,0)[lb]{\smash{{\SetFigFont{12}{14.4}{\rmdefault}{\mddefault}{\updefault}{\color[rgb]{0,0,0}$a$}%
}}}}
\put(4681,-1591){\makebox(0,0)[lb]{\smash{{\SetFigFont{12}{14.4}{\rmdefault}{\mddefault}{\updefault}{\color[rgb]{0,0,0}$b$}%
}}}}
\put(2971,-151){\makebox(0,0)[lb]{\smash{{\SetFigFont{12}{14.4}{\rmdefault}{\mddefault}{\updefault}{\color[rgb]{0,0,0}$Q_n$}%
}}}}
\end{picture}%